\documentclass[12pt]{amsart}

\usepackage{amsfonts,amsmath,amssymb}
\title[Isomorphism classes of type $A_n$]
{Isomorphism classes and automorphisms of finite Hopf algebras of type $A_n$}
\author{Nicol\'{a}s Andruskiewitsch}
\address{Facultad de Matem\'{a}tica, Astronom\'{i}a y f\'{i}sica\\
Universidad Nacional de C\'{o}rdoba \\ CIEM - CONICET,
(5000) Ciudad Universitaria \\
C\'{o}rdoba \\Argentina}
\email{andrus@mate.uncor.edu}
\author{Hans-J\"urgen Schneider}
\address{Mathematisches Institut, Universit\"at M\"unchen, Theresienstr. 39, D-80333 Munich,
Germany}
\email{Hans-Juergen.Schneider@mathematik.uni-muenchen.de}
\thanks{Results of this paper were obtained during a visit of
H.-J. S. at the University of C\'ordoba, partially supported
through a grant of CONICET. The work of N. A. was partially
supported by CONICET, Fundaci\' on Antorchas, Agencia C\'ordoba
Ciencia, ANPCyT, Secyt (UNC) and TWAS (Trieste)}

\newtheorem{Lem}{Lemma}[section]

\newtheorem{Cor}[Lem]{Corollary}
\newtheorem{Thm}[Lem]{Theorem}

\theoremstyle{definition}
\newtheorem{Def}[Lem]{Definition}

\newcommand{\epf}{$\Box$}
\newcommand{\pf}{\medskip\noindent{\sc Proof: }}

\newcommand\Aut{\operatorname{Aut}}

\newcommand\id{{\operatorname{id}}}

\newcommand\ad{{\operatorname{ad}}}

\newcommand\gr{{\operatorname{gr}}}

\newcommand\Isom{{\operatorname{Isom}}}

\renewcommand\o{\otimes}

\newcommand\YDG{{}^{\Gamma}_{\Gamma}\mathcal{YD}}

\newcommand\D{\mathcal{D}}

\newcommand\G{\Gamma}

\numberwithin{equation}{section}

\begin{document}

\maketitle

\section{Introduction}

In \cite{AS2} we classified a large class of finite-dimensional
pointed Hopf algebras up to isomorphism. However the following
problem was left open for Hopf algebras of type $A,D$ or $E_6$,
that is whose Cartan matrix is connected and allows a non-trivial
automorphism of the corresponding Dynkin diagram. In this case we
described the isomorphisms between two such Hopf algebras with the
same Cartan matrix only implicitly. The problem is whether it is
possible to compute the isomorphisms in terms of the defining
families of parameters.

In the present paper we solve this problem for type $A$. To our
surprise there are closed formulas for these isomorphisms. They
are based on an action of the non-trivial automorphism $\sigma$ of
the Dynkin diagram on the parameter spaces of the Hopf algebras of
type $A$.

\medskip

The Hopf algebras $u(\D,\mu)$ of type $A_n$ can be defined as
follows. For more details and references to the literature we
refer to our survey paper \cite{AS1}. Let $n \geq 2$ and
$(a_{ij})_{1 \leq i,j \leq n}$ the Cartan matrix of type $A_n$ in
the form
\begin{equation}\label{standard}
a_{ij} =
\begin{cases}
2, & \text{ if } i=j,\\
-1, & \text{ if }|i-j|=1,\\
0, &  \text{ if }|i-j| >1.
\end{cases}
\end{equation}
Let $\G$ be a finite abelian group, $g_i \in \G$ and $\chi_i$
characters of $\G$ for all $1 \leq i \leq n$. Define $q_{ij} =
\chi_j(g_i),1 \leq i,j \leq n.$ Then
$$\D = \D(\G,(g_i)_{1 \leq i \leq n}, (\chi_i)_{1 \leq i \leq n}, (a_{ij})_{1 \leq i,j \leq n})$$
is a {\em datum of Cartan type} if there is a root of unity $q$ of
order $N >1$ in $k$ such that

\begin{align}
q_{ii} &= q \text{ for all }1 \leq i \leq n, \text{ and }\label{qii}\\
q_{ij}q_{ji} &=
\begin{cases}
q^{-1}, &\text{ if } |i-j|=1,\\
1, & \text{ if } |i-j| >1.
\end{cases}\label{qij}
\end{align}

For simplicity we assume that $N$ is {\em odd}. The case when $N$ is even could be treated in the same way.

Let $\Phi^+$ be the positive roots of the root system of type
$A_n$, and let $k^{\Phi^+}$ be the set of all families $\mu =
(\mu_{ij})_{1 \leq i<j \leq n+1}$ of scalars in $k$. A {\em family
of root vector parameters} for $\D$ is a family $\mu \in
k^{\Phi^+}$ satisfying the following two conditions.
\begin{alignat*}{2}
&\textup{\textup{(R1)}}\quad&&\mu_{ij}=0\text{ for all }1 \leq i<j \leq n+1\text{ with } (g_ig_{i+1} \cdots g_{j-1})^N =1.\\
&\textup{(R2)}&&\mu_{ij}=0\text{ for all }1 \leq i<j \leq n+1\text{ with }(\chi_i \chi_{i+1}\cdots \chi_{j-1})^N \neq 1.
\end{alignat*}
In \eqref{defu} we associate to any family $\mu \in k^{\Phi^+}$
satisfying (R2) a family $(u_{ij}(\mu))_{1 \leq i<j \leq n+1}$ of
elements in the group algebra $k[\G]$. If $\mu$ satisfies (R2) we
can always normalize it such that $\mu$ becomes a family of root
vector parameters without changing the elements $u_{ij}(\mu)$.
This normalization process is discussed in Lemma
\ref{normalization}. The Hopf algebra $u(\D,\mu)$ is generated as
an algebra by the group $\G,$ that is, by generators of $\G$
satisfying the relations of the group, and $x_1,\dots,x_n,$ with
the relations:
\begin{align*}
&\text{({\em Action of the group}) }& &gx_i g^{-1} = \chi_i(g) x_i,  \text{ for all }
i, \text{ and all } g \in \G,&\\
&\text{({\em Serre relations}) }& &\ad_c(x_i)^{1 - a_{ij}}(x_j) = 0, \text{ for all }  i \neq j, &\\
&\text{({\em Root vector relations}) }&& x_{ij}^{N} = u_{ij}(\mu), \text{ for all }1 \leq i<j \leq n+1.
\end{align*}
The coalgebra structure is given by
\begin{align*}
&\Delta (x_i) = g_i \o x_i + x_i \o 1, &&\Delta(g) = g \o g, \text{ for all }1 \leq i \leq \theta,  g \in \G.&
\end{align*}
The Serre relations are the deformed Serre relations where
$$(\ad_c(x_i))(x_{j_1} \cdots x_{j_s}) = x_ix_{j_1} \cdots x_{j_s}-
q_{ij_1} \cdots q_{ij_s}x_{j_1} \cdots x_{j_s}x_i, s \geq 1,$$ is
the braided adjoint action. The root vectors $x_{ij}$ are iterated
braided commutators. They are defined in \eqref{rootvector}.

The non-trivial automorphism of the Dynkin diagram of $A_n$ is the
permutation $\sigma \in \mathbb{S}_n$ defined by $\sigma(i) = n -i
+1$ for all $1 \leq i \leq n$. For each $\D$ we have an action of
$\sigma$ on the parameter spaces by an explicitly defined morphism
of affine algebraic varieties
$$\sigma^{\D} : k^ {\Phi^+} \to k^ {\Phi^+},\mu \mapsto (\sigma_{ij}^{\D}(\mu))_{1 \leq i<j \leq n+1}.$$
The polynomials $\sigma_{ij}^{\D}(\mu)$ of degree $j-i$ are
defined in \eqref{defmu'}. By Theorem \ref{square} they  define an
isomorphism of affine algebraic varieties between the subspaces of
all elements of $k^{\Phi^+}$ satisfying \textup{(R2)} for $\D$
resp. for $\D^{\sigma}$. Here
$$\D^{\sigma} = \D(\G,(g_{\sigma(i)})_{1 \leq i \leq n}, (\chi_{\sigma(i)})_{1 \leq i \leq n},(a_{ij})_{1\leq i,j \leq n}).$$
In Corollary \ref{sigmaiso} we show
$$u(\D^{\sigma},\sigma^{\D}(\mu)) \cong u(\D,\mu)$$
for all $\mu \in k^{\Phi^+}$ satisfying (R2).

Our main result is Theorem \ref{Hopfiso}, where we compute all
Hopf algebra isomorphisms between two Hopf algebras $u(\D',\mu')$
and $u(\D,\mu)$ of type $A_n$. The polynomials $\sigma^{\D}_{ij}$
play an important role in this theorem. The first essential steps
in the proof of Theorem \ref{Hopfiso} is Theorem
\ref{mainreverse}, where we compute the basis representation of
the $N$-th powers of the ``reverse root vectors'' in the usual
PBW-basis formed by the root vectors. The second essential step is
Theorem \ref{mainsystem}, where we prove that the images of the
$N$-th powers of the reverse root vectors in $u(\D,\mu)$ are the
elements $u_{ij}^{\D^{\sigma}}(\sigma^{\D}(\mu))$.

The authors thank the referee for helpful remarks.

\section{Finite Hopf algebras of type $A_n$}\label{sectionAn}

\subsection{Diagrams of type $A_n$, root vectors, and reverse root vectors}

Let $n \geq 2$ and $(a_{ij})_{1 \leq i,j \leq n}$ the Cartan
matrix of type $A_n$ in the form \eqref{standard}. Let
$\mathbb{Z}[I]$ be the free abelian group with basis
$\alpha_1,\dots,\alpha_n$. The Weyl group $W$ of $(a_{ij})$ is the
subgroup of $\Aut(\mathbb{Z}[I])$ generated by the simple
reflections $s_1,\dots,s_n$ defined by $s_i(\alpha_j) = \alpha_j -
a_{ij}\alpha_i$ for all $ 1 \leq i,j \leq n.$ The root system
$\Phi$ of $(a_{ij})$ is defined by $\Phi = \cup_{i=1}^n
W(\alpha_i)$. It has the basis $\alpha_1,\dots,\alpha_n$, and the
positive roots with respect to this basis are the elements
$$\alpha_{ij}=\sum_{l=i}^{j-1} \alpha_l, 1 \leq i <j \leq n+1.$$
Let $w_0$ be the longest element in $W$. We choose the reduced
representation
$$w_0 = s_1s_2 \cdots s_n s_1s_2 \cdots s_{n-1}s_1s_2 \cdots s_{n-2} \cdots s_1$$
of $w_0$ of length $p=\frac{n(n+1)}{2}$. The corresponding convex
ordering of the positive roots
$$\beta_{l} = s_{i_1} \cdots s_{i_{l-1}}(\alpha_{i_l}), 1 \leq l \leq p,$$
is the lexicographic ordering, that is,
$$\alpha_{12}<\alpha_{13} < \cdots < \alpha_{1, n+1} < \alpha_{23}< \cdots \alpha_{2,n+1} < \cdots <\alpha_{n,n+1}.$$
Let $\G$ be a finite abelian group,
$$\D = \D(\G,(g_i)_{1 \leq i \leq n}, (\chi_i)_{1 \leq i \leq n}, (a_{ij})_{1 \leq i,j \leq n})$$
a Cartan datum and define $q_{ij} = \chi_j(g_i),1 \leq i,j \leq
n.$ Then the Cartan condition $q_{ij}q_{ji} = q_{ii}^{a_{ij}},1
\leq i,j \leq n,$ is equivalent to the following: There is a root
of unity $q$ of order $N >1$ in $k$ such that \eqref{qii} and
\eqref{qij} hold.

For simplicity we assume that $N$ is {\em odd}; as we said, this
is not essential and the case when $N$ is even could be treated in
the same way.

Let $V \in\, \YDG$ with basis $x_i \in V_{g_i}^{\chi_i}$, $1 \leq
i \leq n$; that is $g \cdot x_i= \chi_i(g) x_i$  for all $g \in
\Gamma$, and $\delta(x_i) = g_i \otimes x_i$.  Then
$$R=R(\D)= k\langle x_1,\dots,x_n \mid \ad_c(x_i)^{1-a_{ij}}(x_j) =0, \, \forall  1 \leq i,j \leq n, i \neq j \rangle$$
is a Hopf algebra in the braided category $\YDG$.
For $x,y \in R$ we define the braided commutator
$$[x,y]_c = xy - \mu c(x \otimes y),$$
where $c$ denotes the braiding and $\mu$ the multiplication map of
$R$. As in \cite[(6-7) and (6-8)]{AS1} we define root vectors
$x_{ij}, 1 \leq i < j \leq n+1,$ in $R$ inductively by
\begin{align}
x_{i,i+1} &= x_i \text{ for all } 1 \leq i \leq n,\\
x_{ij}&=[x_{i,i+1},x_{i+1,j}]_c \text{ for all } 1 \leq i < j \leq n+1, j-i>1.\\
\intertext{Then}
x_{ij} &= [x_{il},x_{lj}]_c \text{ for all } 1 \leq i < l < j \leq n+1,\label{rootvector}
\end{align}
and the root vectors $x_{ij}$ in the lexicographic order define a
PBW-basis of $R$ \cite[Theorem (6.13)]{AS1}.

In addition we define inductively reverse root vectors $x_{ji}, 1 \leq i<j \leq n+1$, in $R$ by
\begin{align}
x_{i+1,i} &= x_i \text{ for all } 1 \leq i \leq n,\\
x_{ji}&=[x_{j,j-1},x_{j-1,i}]_c \text{ for all } 1 \leq i < j \leq n+1, j-i>1.\label{root}\\
\intertext{Again it follows that}
x_{ji} &= [x_{jl},x_{li}]_c \text{ for all } 1 \leq i < l < j \leq n+1.\label{root1}
\end{align}
Thus for all $1 \leq i<j \leq n+1,$ $x_{ij}$ is any bracketing of
the elements $x_i,x_{i+1},\dots,x_{j-1}$ in this order, and
$x_{ji}$ is any bracketing of the reverse sequence
$x_{j-1},x_{j-2},\dots,x_i$.

\subsection{Root vector parameters  and normalization}

For any positive root we define elements in the group and characters of the group by
\begin{equation}\label{defgij}
g_{ij}=\prod_{i \leq l <j}g_l, \chi_{ij}=\prod_{i \leq l <j}\chi_l, 1 \leq i<j \leq n+1.
\end{equation}
A {\em family of root vector parameters} for $\D$ is a family $\mu
= (\mu_{ij})_{1 \leq i<j \leq n+1}$ of scalars $\mu_{ij} \in k$
satisfying the following two conditions.
\begin{alignat*}{2}
&\textup{\textup{(R1)}}\quad&&\mu_{ij}=0\text{ for all }1 \leq i<j \leq n+1\text{ with } g_{ij}^N =1.\\
&\textup{(R2)}&&\mu_{ij}=0\text{ for all }1 \leq i<j \leq n+1\text{ with }\chi_{ij}^N \neq 1.
\end{alignat*}
For all $ 1 \leq i<j \leq n+1$ let
$$I_{ij} = \{ (i_1,\dots,i_r) \mid r \geq 2, i=i_1 < i_2 < \cdots < i_r = j \}.$$
We denote the set of all families $\mu = (\mu_{ij})_{1 \leq i<j
\leq n+1}$ of elements in $k$ by $k^{\Phi^+}$.

For any $\mu \in k^{\Phi^+}$  we define for all $1 \leq i<j \leq n+1$ scalars
\begin{align}
\mu(i_1,\dots,i_r)&=\mu_{i_1i_2} \cdots \mu_{i_{r-1}i_{r}} \text{ for all } (i_1,\dots,i_r) \in I_{ij},\label{longmu}\\
\intertext{and elements in the group algebra} u^{\D}_{ij}(\mu) &=
\sum_{(i_1,\dots,i_r) \in I_{ij}} (q-1)^{N(r-1)}
\mu(i_1,\dots,i_r) (1-g_{i_{r-1}i_r}^N).\label{defu}
\end{align}
Thus \begin{align*}u^{\D}_{ij}(\mu)&= \mu_{ij}(1-g_{ij}^N) \\ &+
\sum_{i<p<j}\big(\sum_{(i_1,\dots,i_r) \in I_{ip}} (q-1)^{N(r-1)}
\mu(i_1,\dots,i_r)\big) \mu_{pj}(1-g_{pj}^N).\end{align*} We will
write $u_{ij}(\mu)=u^{\D}_{ij}(\mu)$ when $\D$ is fixed. Recall
that $q = \chi_i(g_i)$ for all $1 \leq i \leq n$ also depends on
$\D$. It is easy to see that the family $(u_{ij}(\mu))_{1 \leq i<j
\leq n+1}$ of elements in the group algebra can inductively be
defined by
\begin{equation}\label{inductive}
u_{ij}(\mu) =\mu_{ij}(1-g_{ij}^N) + \sum_{i<p<j}(q-1)^N \mu_{ip}
u_{pj}(\mu), 1 \leq i<j \leq n+1.
\end{equation}
This definition agrees with the inductive definition in \cite[Theorem 6.18]{AS1}
when $\mu$ is a family of root vector parameters for $\D$. There
we defined
\begin{align*}
C_{kl}^j &= (1-q^{-1})^N \chi_{kl}(g_{lj})^{\frac{N(N-1)}{2}},i \leq k<l \leq j,\\
u_{ij}(\mu) &=\mu_{ij}(1-g_{ij}^N) + \sum_{i<p<j}C_{ip}^j \;\mu_{ip}
u_{pj}(\mu), 1 \leq i<j \leq n+1.
\end{align*}
Since $N$ is odd it follows from \textup{(R2)} that $C_{kl}^j
\mu_{kl} = (q-1)^N \mu_{kl}$ for all $i \leq k<l \leq j$. Thus
both definitions do agree.

\begin{Lem}\label{lemmaroot3}
Let $\mu \in k^{\Phi^+}$.
\begin{enumerate}
\item Suppose $\mu$ satisfies \textup{(R2)}. Then
\begin{alignat*}{2}
&\textup{(R3)}\quad &&u_{ij}(\mu)=0\text{ for all }1 \leq i<j \leq n+1 \text{ with } \chi_{ij}^N \neq 1.
\end{alignat*}
\item Suppose $\mu$ satisfies \textup{\textup{(R1)}} and \textup{(R3)}.
Then $\mu$ satisfies \textup{(R2)}, that is, $\mu$ is a family of
root vector parameters for $\D$.
\end{enumerate}
\end{Lem}
\pf This follows by induction on $j-i$ from \eqref{inductive}
since for all $i<p<j$ the inequality $\chi_{ij}^N \neq 1$ implies
that $\chi_{ip}^N \neq 1$ or $\chi_{pj}^N \neq 1$. \epf

By \cite[Theorem 6.18]{AS1} the families $(u_{ij}(\mu))_{1 \leq
i<j \leq n+1}$ are exactly the solutions of the equations
\begin{equation}\label{reason}
\Delta(u_{ij})=
u_{ij} \otimes 1 + g_{ij}^N \otimes u_{ij}
+ \sum_{i<p<j}(q-1)^N u_{ip}g_{pj}^N \otimes u_{pj}
\end{equation}
in $k[\G] \otimes k[\G]$ for all $1 \leq i<j \leq n+1$. This
characterization of the $u_{ij}(\mu)$ is used to prove the next
lemma. It shows how to ``normalize'' an arbitrary sequence $\mu$
so that \textup{\textup{(R1)}} is satisfied.

\begin{Lem}\label{normalization}
Let $\mu \in k^{\Phi^+}$. Then there is exactly one family $\mu'
\in k^{\Phi^+}$ satisfying \textup{\textup{(R1)}} such that
\begin{equation}\label{defnormal}
u_{ij}(\mu) = u_{ij}(\mu') \text{ for all }1 \leq i<j \leq n+1.
\end{equation}
If $\mu$ satisfies \textup{(R2)} then $\mu'$ is a family of root vector parameters for $\D$.
\end{Lem}
\pf
Let $u_{ij}=u_{ij}(\mu)$ for all $ 1 \leq i<j \leq n+1$.
We define the elements $\mu_{ij}'$ by induction on $j-i$.
Let $j=i+1$. Let
$$\mu_{i,i+1}' = \begin{cases}
\mu_{i,i+1}, &\text{ if } g_{i,i+1}^N \neq1,\\
0,&\text{ if } g_{i,i+1}^N =1.
\end{cases}$$
Then \eqref{defnormal} for $(i,i+1)$ holds since $u_{i,i+1}(\mu)=\mu_{i,i+1}(1 - g_{i,i+1}^N)$.

Let $k>1$. Suppose we have already defined $\mu_{ij}'$ whenever
$j-i \leq k-1$ such that \eqref{defnormal} holds if $j-i \leq
k-1$. Let $1 \leq i<j \leq n+1$ and assume that $j-i=k$. If
$g_{ij}^N =1$, then we define $\mu_{ij}'=0$. If $g_{ij}^N \neq 1$,
we define $\mu_{ij}' \in k$ to be the unique scalar satisfying
\begin{equation}\label{proofnormal}
u_{ij} = \mu_{ij}'(1-g_{ij}^N) + \sum_{i<p<j} (q-1)^N \mu_{ip}'u_{pj}.
\end{equation}
The existence of the scalar $\mu_{ij}'$ follows from the argument
in the induction step of the proof of \cite[Theorem 6.18]{AS1}.

Thus we have shown the existence of the family $\mu'$. Uniqueness
follows easily by induction on $j-i$ from \eqref{inductive}.

Suppose that $\mu$ satisfies \textup{(R2)}. Then $\mu$ satisfies
\textup{(R3)} by Lemma \ref{lemmaroot3} (1), and $\mu'$ satisfies
\textup{(R2)} by Lemma \ref{lemmaroot3} (2). \epf

For any $\mu \in k^{\Phi^+}$ we define
\begin{equation}\label{defTheta}
\nu^{\D}(\mu)=\mu',\quad \nu_{ij}^{\D}(\mu)= \mu_{ij}' \text{ for
all }1 \leq i<j \leq +1,
\end{equation}
where $\mu'$ is the family constructed from $\mu$ in Lemma
\ref{normalization}. We call $\nu^{\D}(\mu)$ the {\em
normalization} of $\mu$.

The elements $\mu_{ij}' = \nu^{\D}_{ij}(\mu)$ can be computed
inductively. Let $1 \leq i<j \leq n+1$ and assume that $g_{ij}^N
\neq 1$. Then we have $u_{ij}(\mu) = u_{ij}(\mu')$. We replace
$u_{ij}(\mu)$ and $u_{ij}(\mu')$ by the right hand sides of
\eqref{defu} and collect all terms with coefficient
$(1-g_{ij}^N)$. This gives the equality
\begin{align}
&\mu_{ij}
 + \sum_{\substack{i<p<j\\g_{ip}^N=1}}\sum_{(i_1,\dots,i_r) \in I_{ip}} (q-1)^{N(r-1)} \mu(i_1,\dots,i_r) \mu_{pj}\label{construction}\\
={}&\mu_{ij}'
 + \sum_{\substack{i<p<j\\g_{ip}^N=1}}\sum_{(i_1,\dots,i_r) \in I_{ip}} (q-1)^{N(r-1)} \mu'(i_1,\dots,i_r) \mu'_{pj},\notag
\end{align}
where $\mu'(i_1,\dots,i_r)$ is defined in \eqref{longmu} for
$\mu'$. Thus $\mu_{ij}'$ is a function of $\mu$ and $\mu_{ab}',
b-a<j-i$. In particular, we see that $\nu^{\D}_{ij}$ is a
polynomial in the variables $(\mu_{ij})$ with coefficients in
$\mathbb{Z}[q]$. The polynomial $\nu^{\D}_{ij}$ depends on $\D$,
more precisely on $q=\chi_i(g_i)$ and on the numbers
\begin{alignat*}{2}
d_{ab}&=\begin{cases}
1,& \text{ if }g_{ab}^N\neq 1,
\\0,& \text{ if } g_{ab}^N= 1.
\end{cases}, \qquad b-a \leq j-i.\\
\intertext{For example assume that $g_{i,i+2}^N \neq 1$. Then}
\mu_{i,i+2}' &= \begin{cases}
\mu_{i,i+2},&\text{ if } g_i^N\neq1,\\
\mu_{i,i+2} + \mu(i,i+1,i+2),& \text{ if }g_i^N=1.
\end{cases}
\end{alignat*}

\subsection{The Hopf algebras $u(\D,\mu)$ and isomorphisms}\label{subsectioniso}

As in \cite{AS1,AS2} we define for any family of root vector
parameters $\mu$ for $\D$ a finite-dimensional Hopf algebra by
\begin{equation}\label{defHopf}
u(\D,\mu) = R(\D)\# k[\G]/(x_{ij}^N - u_{ij}(\mu) \mid 1 \leq i<j \leq n+1).
\end{equation}
We can extend this definition to all $\mu \in k^{\Phi^+}$
satisfying \textup{(R2)} since by Lemma \ref{normalization} then
$\nu^{\D}(\mu)$ is a family of root vector parameters, and
$u_{ij}(\mu) = u_{ij}(\nu^{\D}(\mu))$ for all $1 \leq i<j \leq
n+1$. By \cite[Theorem 6.25]{AS1} a finite-dimensional pointed
Hopf algebra $A$ is of the form $A \cong u(\D,\mu)$ if and only if
$\gr(A) \cong u(\D,0) \# k[\G],$ where $\gr(A)$ is the graded Hopf
algebra associated to the coradical filtration of $A$.

Let $\rho \in \mathbb{S}_n$ be a {\em diagram automorphism} of
$(a_{ij})$, that is, $$a_{ij}=a_{\rho(i)\rho(j)} \text{ for all }1
\leq i,j \leq n.$$ Then $ \rho = \id$ or $\rho = \sigma$, where
\begin{equation}\label{defsigma}
\sigma(i) = n-i+1 \text{ for all }1 \leq i \leq n.
\end{equation}
As in \cite[Theorem 7.5]{AS2} let
$$\D^{\rho} = \D(\G,(g^{\rho}_i)_{1 \leq i \leq n}, (\chi^{\rho}_i)_{1 \leq i \leq n},(a_{ij})_{1\leq i,j \leq n})$$
be the Cartan datum with $g^{\rho}_i = g_{\rho(i)}, \chi^{\rho}_i
= \chi_{\rho(i)}$ for all $1 \leq i \leq n$. Let $V^{\rho} \in
{\YDG}$ with basis $x_i^{\rho} \in
(V^{\rho})_{g_{\rho(i)}}^{\chi_{\rho(i)}}$ for all $1 \leq i \leq
n$. Then
$$F^{\rho} : R(\D^{\rho}) \to R(\D), \quad x_i^{\rho} \mapsto x_{\rho(i)} \text{ for all } 1 \leq i \leq n,$$
defines an isomorphism of braided Hopf algebras in $\YDG$. For all
$1 \leq i < j \leq n+1$, we denote the  root vector of
$\alpha_{ij}$ in $R(\D^{\rho})$ by $x_{ij}^{\rho}$.

Let $\D' = \D(\Gamma', (g'_i)_{1 \leq i \leq n}, (\chi'_i)_{1 \leq
i \leq n}, (a_{ij})_{1\leq i,j \leq n})$ be another Cartan datum
with finite abelian group $\G'$ and the same Cartan matrix of type
$A_n$ as $\D$. Let $\varphi : \G' \to \G$ be a group isomorphism,
$\rho \in \mathbb S_{n}$  a diagram automorphism of $(a_{ij})$ and
$s=(s_i)_{1 \leq i \leq n}$ a family of non-zero elements in $k$.
Let
\begin{equation}\label{defs}
s_{ij} =\prod_{i \leq l <j}s_l \text{ for all } 1 \leq i<j \leq n+1.
\end{equation}
Let $\pi : R(\D)\# k[\G] \to u(\D,\mu)$ be the canonical projection.

The triple $(\varphi, \rho, (s_i))$ is called an {\em isomorphism}
from $(\D',\mu')$ to $(\D,\mu)$ if the following conditions are
satisfied:
\begin{align}
\varphi(g'_i)&=g_{\rho(i)}, \chi'_i=\chi_{\rho(i)}\varphi\text{ for all } 1 \leq i \leq n.\label{I1}\\
\varphi(u^{\D'}_{ij}(\mu')) &= s_{ij}^N \pi(F^{\rho}(x_{ij}^{\rho})^N) \text{ for all }1 \leq i<j \leq n+1.\label{I3}
\end{align}
Let $\Isom(\D',\mu'),(\D,\mu)$ be the set of all isomorphisms from
$(\D',\mu')$ to $(\D,\mu)$. For Hopf algebras $A',A$ we denote by
$\Isom(A',A)$ the set of all Hopf algebra isomorphisms from $A'$
to $A$. Then
\begin{Thm}\cite[Theorem 7.2]{AS2}\label{cite}
The map
$$\Isom((\D',\mu'),(\D,\mu)) \to \Isom(u(\D',\mu'),u(\D,\mu))$$
given by $(\varphi, \rho, (s_i)) \mapsto F$, where $F(x'_i) = s_i
x_{\rho(i)}$ and $F(g') = \varphi(g')$ for all $1 \leq i \leq
\theta$ and $g' \in \G'$, is bijective.
\end{Thm}

The main result in this paper is the explicit computation of the
set $\Isom((\D',\mu'),(\D,\mu))$. In Section \ref{sectioninverse}
we first compute the elements $F^{\sigma}(x_{ij}^{\sigma})^N$ in
terms suitable for our purpose. The next lemma shows that these
elements are reverse root vectors. This lemma also allows to
derive \eqref{root1} from \eqref{rootvector}.
\begin{Lem}\label{Lemmaiso1}
For all $1 \leq i<j \leq n+1$,
\begin{equation*}
F^{\sigma}(x_{ij}^{\sigma}) = x_{n-i+2,n-j+2}.
\end{equation*}
\end{Lem}
\pf This follows by induction on $j-i$. Suppose that $j=i+1$. Then
$x_{ij} = x_i$ and $F^{\sigma}(x_{ij}^{\sigma})= x_{\sigma(i)} =
x_{n-i+1} = x_{n-i+2,n-j+2}$.

If $j-i\geq 2$, then
\begin{align*}
F^{\sigma}(x_{ij}^{\sigma})&= F^{\sigma}([x_i^{\sigma},x_{i+1,j}^{\sigma}]_{c^{\sigma}})&&\\
&=[x_{\sigma(i)},F^{\sigma}( x_{i+1,j}^{\sigma})]_c&&\\
&=[x_{n-i+1},x_{n-i+1,n-j+2}]_c&&\text{ (by induction)}\\
&=x_{n-i+2,n-j+2}&& \text{ (by \eqref{root})}.
\end{align*}
\epf

\section{The reverse root vectors}\label{sectioninverse}

In the next theorem we compute the basis representation of the
$N$-th powers of the reverse root vectors in the standard
PBW-basis.

As in the last section we  fix a diagram $\D$ of Cartan type $A_n$
and let $R=R(\D)$. For all $1 \leq i <j \leq n+1$ we define
\begin{align}
\tau_{ij} &=\prod_{i\leq k<l <j} q_{lk}^N.\label{tau}\\
\tau(i_1,\dots,i_r)&=\tau_{i_1i_2}\tau_{i_2i_3} \cdots
\tau_{i_{r-1}i_r},\text{ for all } (i_1,\dots,i_r) \in
I_{ij}.\label{longtau}
\end{align}
Note that $\tau_{ij} =  \prod_{i<l<j} \chi_{il}^N(g_l)$. We write
$\tau^{\D}_{ij}$ instead of $\tau_{ij}$ if we want to emphasize
the datum $\D$.

\begin{Thm}\label{mainreverse}
Assume that $1 \leq i<j \leq n+1$. For all $(i_1,\dots,i_r) \in I_{ij}$ define
\begin{align}
t(i_1,\dots,i_r) &= (-1)^{j-i-r+1}(q-1)^{N(r-2)}\tau(i_1,\dots,i_r)^{-\frac{N-1}{2}}\tau_{ij}^{\frac{N+1}{2}}.\notag\\
\intertext{Then}
x_{ji}^N &= \sum_{(i_1, \dots,i_r) \in I_{ij}} t(i_1,\dots,i_r) x_{i_1i_2}^N x_{i_2i_3}^N \cdots x_{i_{r-1}i_r}^N.\label{claim}
\end{align}
\end{Thm}

The proof of Theorem \ref{mainreverse} will be done after Lemma
\ref{computet1}.

To compute the coefficients $t(i_1,\dots,i_r)$ we first change the
notation using characteristic functions. We can assume that $j-2
\geq i$ since $x_{i+1,i}=x_i$. For natural numbers $k<l$ let
$[k,l]=\{k,k+1,\dots,l\}$.

Let $E_{ij}$ be the set of all functions $e : [i,j-2] \to
\mathbb{N}$ with values in $\{0,1\}$. We consider the bijection
$$\Omega : I_{ij} \to E_{ij}$$
given for all $(i_1,\dots,i_r) \in I_{ij}$ and $l \in [i,j-2]$ by
$$\Omega(i_1,\dots,i_r)(l) = \begin{cases}
1,& \text{ if } l  \in \{i_{2} -1,\dots,i_{r-1} -1\},\\
0,& \text{ otherwise}.
\end{cases}$$
For any $e \in E_{ij}$ define
\begin{align}
|e| &= \#\{l \mid i \leq l \leq j-2, e(l)=0\}.\notag\\
\intertext{If $e=\Omega(i_1,\dots,i_r)$ then}
|e| &=j-i-r+1.\label{|e|}
\end{align}
The constant function in $E_{ij}$ with value 1 resp. 0 will be denoted by $(1)$ resp. $(0)$.
Thus $(1)=\Omega(i,i+1,\dots,j)$ and $(0) = \Omega(i,j)$.

\medbreak For $e,f \in E_{ij}$ we write $e \leq f$ if for all $i
\leq l \leq j-2,e(l)=0$ implies $f(l)=0$.

\begin{Lem}\label{ef}
Let $1 \leq i<j \leq n+1, j-i \geq 2$, and $f \in E_{ij}, (1) \neq f$. Then
\begin{align}
\sum_{e \in E_{ij},e \leq f} (-1)^{|e|} &=0,\label{ef1}\\
\sum_{(i_1,\dots,i_r) \in I_{ij}} (-1)^r &=0.\label{ef2}
\end{align}
\end{Lem}
\pf Since $f \neq (1)$ we can choose an index $l$ with $f(l) =0$.
Then $\{e \in E_{ij} \mid e \leq f\}$ is the disjoint union of
elements $e$ with $e(l)=0$ and with  $e(l)=1$, and \eqref{ef1} is
obvious. To prove \eqref{ef2} we consider the case of \eqref{ef1}
with $f=(0)$. Then $\sum_{e \in E_{ij}} (-1)^{|e|} =0$, and
\eqref{ef2} follows from the bijection $\Omega$ and \eqref{|e|}.
\epf

For all $e=\Omega(i_1,\dots,i_r) \in E_{ij}$ let
\begin{align*}
\tau_e&= (q-1)^{N|e|} (\tau_{i_1i_2} \tau_{i_2i_3} \cdots \tau_{i_{r-1}i_r})^{\frac{N-1}{2}},\\
t_e &= (-1)^{|e|} \tau_e^{-1} (q-1)^{N(j-i-1)}\tau_{ij}^{\frac{N+1}{2}},\\
x_e^N &= x_{i_1i_2}^N x_{i_2i_3}^N \cdots x_{i_{r-1}i_r}^N.\\
\end{align*}
Note that
$t_e = t(i_1,\dots,i_r)\text{ if } e=\Omega(i_1,\dots,i_r)$.
This follows from the definitions using \eqref{|e|}.
Hence \eqref{claim} in Theorem \ref{mainreverse} can be restated as
\begin{equation}\label{claim1}
x_{ji}^N = \sum_{e \in E_{ij}} t_ex_e^N.
\end{equation}
The idea of the proof of Theorem \ref{mainreverse} is to project
$R$ onto skew-polynomial rings $R_e$, one for each $e \in E_{ij}$.
Before we begin the proof we establish some technical results on
these projections.
\begin{Def}\label{DefR_a}
For any $e \in E_{ij}$  let $R_e$ be the algebra generated by $x_i, x_{i+1},\dots,x_{j-1}$ with relations
\begin{align}
x_lx_{l+1} - q_{l,l+1}x_{l+1}x_l&= 0,\text{ if } e(l)=1,i \leq l \leq j-2,\label{a=1}\\
x_{l+1}x_l - q_{l+1,l}x_lx_{l+1} &=0,\text{ if } e(l)=0,i \leq l \leq j-2,\label{a=0}\\
x_kx_l - q_{kl} x_lx_k &=0, \text{ if } i \leq k,l \leq j-1, |k-l| \geq 2.\label{skew}
\end{align}
\end{Def}

\begin{Lem}\label{projection}
For any $f \in E_{ij}$, the natural projection
$$\pi_f : R \to R_f, \pi_f(x_l)=
\begin{cases}
x_l,& \text{ if } i \leq l \leq j-1,\\
0,& \text{ otherwise },
\end{cases}, 1 \leq l \leq n,$$
is a well-defined algebra map, and
for all $i \leq u < v \leq j, v-u \geq 2,$
\begin{align}
\pi_f(x_{uv}) =0, \text{ if } f(l)=1 \text{ for some } u \leq l < l+2 \leq v,\label{f(root)}\\
\pi_f(x_{vu}) =0, \text{ if } f(l)=0 \text{ for some } u \leq l < l+2 \leq v.\label{f(root')}
\end{align}
\end{Lem}
\pf
The Serre relations can be reformulated according to the following identities
\begin{align}
-q_{l+1,l} \ad_c(x_l)^2(x_{l+1})&= x_l[x_{l+1},x_l]_c - q_{l,l+1}[x_{l+1},x_l]_cx_l,\\
-q_{l,l+1}\ad_c(x_{l+1})^2(x_{l})&= x_{l+1}[x_l,x_{l+1}]_c - q_{l+1,l} [x_l,x_{l+1}]_cx_{l+1},
\end{align}
in the free algebra $k\langle x_1,\dots,x_n\rangle$ for all $1
\leq l \leq n-1$. Hence both Serre relations
$\ad_c(x_l)^2(x_{l+1}) =0$ and $\ad_c(x_{l+1})^2(x_{l}) =0$ hold
in $R_f$ for all $i \leq l <j-1$, since  $[x_l,x_{l+1}]_c=0$ by
\eqref{a=1} if $f(l)=1$, and $[x_{l+1},x_l]_c=0$ by \eqref{a=0} if
$f(l)=0$. Thus $\pi_f$ is well-defined.

To prove \eqref{f(root)} note that by \eqref{rootvector}
$$x_{uv} = \begin{cases}
[x_{ul},x_{l,l+2}]_c, &\text{ if }u < l < l+2 =v,\\
[x_{ul},[x_{l,l+2},x_{l+2,v}]_c]_c , & \text{ if } u < l < l+2 <v,
\end{cases}$$
and that $\pi_f(x_{l,l+2})=0$ by definition of $\pi_f$. In the
same way \eqref{f(root')} follows from \eqref{root1}. \epf

\medbreak We note the following obvious rule in skew polynomial
rings.

\begin{Lem}\label{skewrule}
Let $x_1,\dots,x_m$ be elements in an algebra such that
$$x_lx_k=p_{lk} x_kx_l \text{ for all }k<l,$$
where $p_{kl} \in k$ for all $k<l$.
Then for any natural number $N$,
$$(x_1 \cdots x_m)^N = \prod_{k<l}p_{lk}^{\frac{N(N-1)}{2}} x_1^N \cdots x_m^N.$$
\end{Lem}
\epf

\begin{Lem}\label{computationkl} Let $f \in E_{ij}$ and $i \leq k <l \leq j$ with $k \leq l-2$.
Suppose that $f(k) = f(k+1)=\cdots=f(l-2) = 0$. Then
$$\pi_f(x_{kl}^N) =(q-1)^{N(l-k-1)}\tau_{kl}^{\frac{N-1}{2}}x_k^Nx_{k+1}^N \cdots x_{l-1}^N.$$
\end{Lem}
\pf
We first prove by induction on $l-k$ that
\begin{equation}\label{kl1}
\pi_f(x_{kl})=(1 - q^{-1})^{l-k-1} x_k x_{k+1} \cdots x_{l-1}.
\end{equation}
Suppose that $l=k+2$. Then
\begin{align*}
\pi_f(x_{k,k+2}) &= x_kx_{k+1} - q_{k,k+1}x_{k+1}x_k &&\\
&=x_kx_{k+1} - q_{k,k+1}q_{k+1,k}x_kx_{k+1}&& \text{ (by \eqref{a=0} and  $f(k)=0$)}\\
&= (1-q^{-1})x_kx_{k+1}&& \text{ (by \eqref{qij})}.
\end{align*}
This proves \eqref{kl1} for $l-k=2$.

For the induction step let $l-k>2$. We obtain by induction
\begin{align*}
\pi_f(x_{kl})&=[x_k,x_{k+1,l}]_c \\
&=[x_k,(1-q^{-1})^{l-k-2} x_{k+1}x_{k+2} \cdots x_{l-1}]_c \\
&=(1-q^{-1})^{l-k-2}(x_kx_{k+1}\cdots x_{l-1}\\
 &\phantom{=(1-q^{-1})^{l-k-2}}- q_{k,k+1}q_{k,k+2} \cdots q_{k,l-1} x_{k+1}x_{k+2} \cdots x_{l-1}x_k).
\end{align*}
Since
$x_{k+1}x_{k+2} \cdots x_{l-1}x_k=q_{l-1,k}q_{l-2,k} \cdots q_{k+1,k}x_kx_{k+1}\cdots x_{l-1}$
in $R_f$ by $f(k)=0$, \eqref{a=0} and \eqref{skew}, and
$$q_{k,k+1}q_{k,k+2} \cdots q_{k,l-1}q_{l-1,k}q_{l-2,k} \cdots q_{k+1,k}=q^{-1}$$
by \eqref{qij}, equation \eqref{kl1} for $k,l$ follows.

\medskip
Since $f(k) = f(k+1)=\cdots=f(l-2) = 0$, \eqref{a=0} and
\eqref{skew} imply by Lemma \ref{skewrule} that
$$(x_kx_{k+1} \cdots x_{l-1})^N = \prod_{k \leq \mu < \nu <l} q_{\nu,\mu}^{\frac{N(N-1)}{2}} x_k^N x_{k+1}^N \cdots x_{l-1}^N.$$
Hence Lemma \ref{computationkl} follows from \eqref{kl1} by taking $N$-th powers.
Note that $(1 -q^{-1})^N = (q-1)^N$ since $q^N=1$.
\epf

\begin{Lem}\label{computepi}
Let $e,f \in E_{ij}$. Then
$$\pi_f(x_e^N)= \begin{cases}
\tau_e\, x_i^N x_{i+1}^N \cdots x_{j-1}^N,& \text{ if } e \leq f,\\
0,& \text{ otherwise}.
\end{cases}$$
\end{Lem}
\pf Let $e=\Omega(i_1,\dots,i_r)$. Suppose that $e \nleq f$, that
is $f(l)=1$ and $e(l)=0$ for some $i \leq l \leq j-2$. Then $l+1
\notin \{i_1,\dots,i_r\}$, since $e(l)=0$. Hence there is an index
$s$ with $i_s \leq l < l+2 \leq i_{s+1}$. Since $f(l)=1$, it
follows by \eqref{f(root)} that $\pi_f(x_{i_s,i_{s+1}}) =0$ , and
thus $\pi_f(x_e^N)=0$.

Now assume $e \leq f$. Since $\pi_f$ is an algebra map, it is
enough to show for all $1 \leq s<r$ that
\begin{equation}\label{computekl} \pi_f(x_{i_s,i_{s+1}}^N)=
(q-1)^{N(i_{s+1}-i_s -1)}\, \tau_{i_s,i_{s+1}}^{\frac{N-1}2}\,
x_{i_s}^N x_{i_s+1}^N\cdots x_{i_{s+1}-1}^N.
\end{equation}
Note that by \eqref{|e|} $\sum_{s=1}^{r-1}(i_{s+1}-i_s -1) = j-i-r+1 = |e|$.

If $i_s +1 = i_{s+1}$, then $x_{i_s,i_{s+1}}=x_{i_s}$ and
\eqref{computekl} is obvious.

If $i_s \leq i_{s+1}-2$, then $e(l)=0$ for all $i_s \leq l \leq
i_{s+1}-2$ by definition of the function $\Omega$. Hence $f(l) =0$
for all $i_s \leq l \leq i_{s+1}-2$ since $e \leq f$, and
\eqref{computekl} follows from Lemma \ref{computationkl}. \epf

\begin{Lem}\label{computet1}
$$\pi_{(1)}(x_{ji}^N) =t_{(1)} x_i^N x_{i+1}^N \cdots x_{j-1}^N.$$
\end{Lem}
\pf
We first prove by induction on $j-i$ that
\begin{equation}\label{computet11}
\pi_{(1)}(x_{ji})= (q-1)^{j-i-1}\left(\prod_{i \leq k<l<j}
q_{lk}\right)\, x_i x_{i+1} \cdots x_{j-1}.
\end{equation}

Suppose that $j=i+2$. Then
\begin{equation}\label{computeR1}
x_{i+1}x_i = q_{i,i+1}^{-1} x_ix_{i+1} \text{ in }R_{(1)}.
\end{equation}
Hence
\begin{align*}
\pi_{(1)}(x_{i+2,i}) &= x_{i+1}x_i - q_{i+1,i}x_ix_{i+1}&&\\
&=(q_{i,i+1}^{-1} - q_{i+1,i})x_ix_{i+1}&&\\
&= (q-1)q_{i+1,i} x_ix_{i+1} &&\text{ (since $q_{i,i+1}q_{i+1,i} = q^{-1}$ by \eqref{qij})}.
\end{align*}
For the induction step let $j-i >2$. Then by induction
\begin{align*}
\pi_{(1)}(x_{ji})&=[x_{j-1},x_{j-1,i}]_c\\
&=[x_{j-1},(q-1)^{j-i-2} \Big(\prod_{i \leq k<l<j-1} q_{lk}\Big)\, x_i x_{i+1} \cdots x_{j-2}]_c\\
&=(q-1)^{j-i-2}\Big(\prod_{i \leq k<l<j-1} q_{lk}\Big)(x_{j-1}x_i x_{i+1} \cdots x_{j-2}\\
&\phantom{=(q-1)^{j-i-2}}-q_{j-1,i}q_{j-1,i+1} \cdots q_{j-1,j-2}x_i x_{i+1} \cdots x_{j-2}x_{j-1}).\\
\end{align*}
Since in $R_{(1)}$
\begin{align}
x_{j-1}x_{j-2} &= q_{j-2,j-1}^{-1} x_{j-2}x_{j-1}\label{ruleR11}\\
&=q_{j-1,j-2}q x_{j-2}x_{j-1},\notag\\
\intertext{and hence}
x_{j-1}x_i x_{i+1} \cdots x_{j-2}&= q_{j-1,i} q_{j-1,i+1} \cdots q_{j-1,j-2}q x_i x_{i+1} \cdots x_{j-1},\notag
\end{align}
it follows that
\begin{align*}
\pi_{(1)}(x_{ji})&=(q-1)^{j-i-2}\prod_{i \leq k<l<j-1} q_{lk} \prod_{k=i}^{j-2} q_{j-1,k}(q-1)x_i x_{i+1} \cdots x_{j-1}\\
&=(q-1)^{j-i-1}\prod_{i \leq k<l<j} q_{lk} x_i x_{i+1} \cdots x_{j-1}.
\end{align*}
This finishes the proof of \eqref{computet11}.

By \eqref{a=1}, \eqref{skew} and Lemma \ref{skewrule}
$$(x_ix_{i+1} \cdots x_{j-1})^N = q^{\frac{N(N-1)}{2}(j-i-1)}
\Big(\prod_{i \leq k<l <j}q_{lk}^{\frac{N(N-1)}{2}}\Big)\,
x_i^Nx_{i+1}^N \cdots x_{j-1}^N.$$ Hence Lemma \ref{computet1}
follows from \eqref{computet11} by taking $N$-th powers. Note that
$q^{\frac{N(N-1)}{2}}=1$ since $N$ is odd by assumption. \epf

\bigbreak We now prove Theorem \ref{mainreverse}.

\pf Since the root vectors $x_{ij}$ in the lexicographic order
define a PBW-basis of $R$, there are uniquely determined
coefficients $\widetilde{t}_e \in k, e \in E_{ij}$, with
\begin{equation}\label{constantswidetilde}
x_{ji}^N = \sum_{e \in E_{ij}} \widetilde{t}_e x_e^N,
\end{equation}
\emph{cf.} \cite[Th. 2.6 (2)]{AS2}-- and compare with \cite[Lemma
6.9]{AS1}. By \eqref{claim1} we have to show that $\widetilde{t}_e
= t_e$ for all $e \in E_{ij}$.

To prove that $\widetilde{t}_{(1)} = t_{(1)}$, we apply
$\pi_{(1)}$ to both sides of \eqref{constantswidetilde}. For all
$(1) \neq e \in E_{ij}$ we see from Lemma \ref{computepi} that
$\pi_{(1)}(x_e^N) =0$, since $e \nleq (1)$. Hence
$\pi_{(1)}(\sum_{e \in E_{ij}} \widetilde{t}_e x_e^N) =
\widetilde{t}_{(1)} x_i^N x_{i+1}^N \cdots x_{j-i}^N$, and
$\widetilde{t}_{(1)} = t_{(1)}$ by Lemma \ref{computet1}.

Let $(1) \neq f \in E_{ij}$. Then $f(l) =0$ for some $i \leq l
\leq j-2$ and $\pi_f(x_{ji})=0$ by \eqref{f(root')}. Hence
applying $\pi_f$ to both sides of \eqref{constantswidetilde} and
using Lemma \ref{computepi} we obtain $0=\sum_{e \leq f}
\widetilde{t}_e \tau_e x_i^N x_{i+1}^N \cdots x_{j-1}^N$, hence
\begin{equation}\label{ttau}
\sum_{e \leq f} \widetilde{t}_e \tau_e=0.
\end{equation}
Note that by definition
\begin{equation*}
t_f = (-1)^{|f|} \tau_f^{-1}t_{(1)}\text{ for all }  f \in E_{ij}.
\end{equation*}
To finish the proof of the theorem we therefore show by induction
on $|f|$ that \begin{equation}\label{formula ttau}
\widetilde{t}_f\tau_f = (-1)^{|f|} t_{(1)} \text{ for all }  f \in
E_{ij}.
\end{equation}
Suppose that $|f|=0$. Then $f=(1)$. Since $\widetilde{t}_{(1)} =
t_{(1)}$ and $\tau_{(1)} =1,$ \eqref{formula ttau} follows for
$f=(1)$.

For the induction step we note that $|e| < |f|$ for all $e \leq f, e \neq f$.
Hence we get by induction from \eqref{ttau} for all $f \neq (1)$
\begin{equation}\label{t3}
\widetilde{t}_f \tau_f = -\sum_{e \leq f, e \neq f} \widetilde{t}_e\tau_e = -\sum_{e \leq f, e \neq f} (-1)^{|e|} t_{(1)}.
\end{equation}
By \eqref{ef1} $\sum_{e \leq f} (-1)^{|e|} =0$ for $f \neq (1)$,
and \eqref{formula ttau} follows from \eqref{t3}. \epf

By the same proof and \eqref{kl1}, \eqref{computet11} but without
taking $N$-th powers we get the basis representation of $x_{ji}$.

\begin{Thm}
Assume that $1 \leq i<j \leq n+1$. Then
$$x_{ji} = (-q)^{j-i-1} \Big(\prod_{i \leq k<l <j} q_{lk}\Big)\,
 \sum_{(i_1,\dots,i_r) \in I_{ij}} (q^{-1}-1)^{r-2} x_{i_1i_2}\cdots x_{i_{r-1}i_r}.$$
\end{Thm}
\epf


\section{The action of the diagram automorphism on root vector parameters}\label{sectionaction}

As in Section \ref{sectionAn} let
\begin{align*}
\D = \D(\G,(g_i)_{1 \leq i \leq n}, (\chi_i)_{1 \leq i \leq n}, (a_{ij})_{1 \leq i,j \leq n})
\end{align*}
be a datum with finite abelian group $\G$ and Cartan matrix
\eqref{standard} of Type $A_n$. Recall that $\sigma$ denotes the
non-trivial diagram automorphism of $(a_{ij})$ given by
\eqref{defsigma}. In this section we will construct for each $\mu
\in k^{\Phi^+}$ satisfying \textup{(R2)} a family
$\sigma^{\D}(\mu)\in k^{\Phi^+}$ satisfying \textup{(R2)} such
that the isomorphism
$$F^{\sigma}\#\id : R(\D^{\sigma})\#k[\G] \to R(\D)\#k[\G]$$
induces an isomorphism $u(\D^{\sigma},\sigma^{\D}(\mu)) \to u(\D,\mu)$.
We begin with a technical lemma to simplify the constants
$\tau(i_1,\dots,i_r)$ in Theorem \ref{mainreverse} when they
appear as factors of certain root vector parameters.

\begin{Lem}\label{simplify}
Let $1 \leq i <j \leq n+1, (l_1,\dots,l_m) \in I_{ij},(k_1,\dots,k_r) \in I_{1m}$. Let $\mu$ be a family of root vector parameters for $\D$. Then
\begin{align*}
\mu(l_1,\dots,l_m)\tau_{ij}=\mu(l_1,\dots,l_m) \tau(l_{k_1},l_{k_2},\dots,l_{k_r}).
\end{align*}
\end{Lem}
\pf The lemma is trivial if one of the factors $\mu_{l_kl_{k+1}}$
of $\mu(l_1,\dots,l_m)$ is zero. Assume that $\mu_{l_kl_{k+1}}
\neq 0$ for all $1 \leq k <m$. Then $(\chi_{l_kl_{k+1}})^N=1$ for
all $1 \leq k <m$. In particular $(\chi_{l_{k_1}l_{k_s}})^N=1$ for
all $2 \leq s \leq r$. Therefore
\begin{align*}
\tau_{ij} &=\prod_{i<l<j} (\chi_{il})^N(g_l)\\
&= \prod_{l_{k_1}<l<l_{k_2}} (\chi_{l_{k_1}l})^N(g_l)\prod_{l_{k_2} \leq l <l_{k_3}}(\chi_{l_{k_1}l})^N(g_l)\cdots \prod_{l_{k_{r-1}} \leq l < l_{k_r}} (\chi_{l_{k_1}l})(g_l)\\
&=\tau_{l_{k_1}l_{k_2}} \cdots \tau_{l_{k_{r-1}}l_{k_r}} = \tau(l_{k_1},l_{k_2},\dots,l_{k_r}),
\end{align*}
since $(\chi_{l_{k_1}l})^N =
(\chi_{l_{k_1}l_{k_s}})^N(\chi_{l_{k_s}l})^N =(\chi_{l_{k_s}l})^N$
for all $2 \leq s <r$ and $l_{k_s} \leq l$. \epf

\bigbreak We also need the following identity in group algebras.

\begin{Lem}\label{group}
Let $G$ be a group, $m\geq2$ and $h_{st}$ elements in $G$ for all $1 \leq s<t \leq m$. Assume
$$h_{rs}h_{s,s+1} = h_{r,s+1} \text{ for all } 1\leq r \leq 2, 1 \leq s < m.$$

Then
\begin{align*}
&\sum_{(k_1,\dots,k_r) \in I_{1m}}(-1)^{r}
(1-h_{k_{r-1}k_r})\notag\\
=&\sum_{(k_1,\dots,k_r) \in I_{1m}} (-1)^r(1-h_{k_{1},k_{1}+1}) \cdots (1-h_{k_{r-1},k_{r-1}+1}).
\end{align*}
\end{Lem}
\pf Let $$S_m=\sum_{(k_1,\dots,k_r) \in I_{1m}}
(-1)^r(1-h_{k_{1},k_{1}+1}) \cdots (1-h_{k_{r-1},k_{r-1}+1}).$$
The Lemma is true for $m=2$. We show by induction on $m>2$ that
\begin{equation}\label{h1m}
S_m=h_{2m}-h_{1m} \text{ if } m>2.
\end{equation}
This holds for $m=3$ since
$$S_3=-(1-h_{12})(1-h_{23}) + 1-h_{12} = (1-h_{12})h_{23} = h_{23}-h_{13}.$$
The induction step follows from
\begin{align*}
S_{m+1} &=\sum_{\substack{(k_1,\dots,k_r) \in I_{1,m+1}\\k_{r-1} =m}} (-1)^r(1-h_{k_{1},k_{1}+1}) \cdots (1-h_{k_{r-1},k_{r-1}+1})\\
&+\sum_{\substack{(k_1,\dots,k_r) \in I_{1,m+1}\\k_{r-1} <m}} (-1)^r(1-h_{k_{1},k_{1}+1}) \cdots (1-h_{k_{r-1},k_{r-1}+1})\\
&=-S_m(1-h_{m,m+1}) +S_m\\
&=S_mh_{m,m+1}.
\end{align*}
On the other hand if $m>2$ then
\begin{align*}
\sum_{(k_1,\dots,k_r) \in I_{1m}}(-1)^{r}
&(1-h_{k_{r-1}k_r})\\
&=1-h_{1m} + \sum_{1<p<m}\sum_{\substack{(k_1,\dots,k_r) \in I_{1m}\\k_{r-1}=p}}(-1)^{r}(1-h_{k_{r-1}k_r})\\
&=1-h_{1m} +  \sum_{1<p<m}\sum_{(k_1,\dots,k_r) \in I_{1p}}(-1)^{r-1}(1-h_{pm})\\
&=1-h_{1m} -(1-h_{2m}) \text{  by \eqref{ef2}}\\
&=h_{2m}-h_{1m}.
\end{align*}
\epf

We introduce the notation
\begin{equation*}
\widetilde{i} = n-i+2 \text{ for all } 1 \leq i \leq n+1
\end{equation*}
for the non-trivial diagram automorphism of $A_{n+1}$.
Note that  the map
\begin{equation}\label{widetilde}
I_{\widetilde{j}\;\widetilde{i}} \to I_{ij},(i_1,\dots,i_r) \mapsto (\widetilde{i_r}, \dots, \widetilde{i_1}),
\end{equation}
is bijective for all $1 \leq i<j \leq n+1$. Recall that
$$\D^{\sigma} = \D(\G,(g^{\sigma}_i)_{1 \leq i \leq n}, (\chi^{\sigma}_i)_{1 \leq i \leq n},(a_{ij})_{1\leq i,j \leq n}),$$
where $g^{\sigma}_i = g_{\sigma(i)}, \chi^{\sigma}_i = \chi_{\sigma(i)}$ for all $1 \leq i \leq n$.
Then
\begin{equation}\label{rule}
g^{\sigma}_{ij} = \prod_{i\leq l<j} g^{\sigma}_l = g_{\widetilde{j}\;\widetilde{i}},
\text{ for all }1 \leq i<j \leq n+1,
\end{equation}
since $g^{\sigma}_{ij} = g^{\sigma}_i g^{\sigma}_{i+1} \cdots
g^{\sigma}_{j-1} = g_{n-i+1}g_{n-i} \cdots
g_{n-j+2}=g_{\widetilde{j}\;\widetilde{i}}.$ For all $ \mu \in
k^{\Phi^+}$ and all $1 \leq i < j \leq n+1$ we define
\begin{align}\label{defmu'}
\sigma^{\D}_{ij}(\mu)&=\tau_{\widetilde{j}\; \widetilde{i}}(-1)^{j-i+1}\sum_{(i_1,\dots,i_r) \in I_{ij}} (q-1)^{N(r-2)} \mu(\widetilde{i_r}, \cdots, \widetilde{i_1}),\\
\sigma^{\D}(\mu) &= (\sigma^{\D}_{ij}(\mu))_{1 \leq i<j \leq n+1}.\notag
\end{align}

Here $\tau_{\widetilde{j}\;
\widetilde{i}}=\tau^{\D}_{\widetilde{j}\; \widetilde{i}}$ and $q$
depend on $\D$, or more precisely on the braiding matrix
$(q_{ij})$ of $\D$. Note that
$q=\chi_i(g_i)=\chi^{\sigma}_i(g^{\sigma}_i)$ for all $1 \leq i
\leq n$.

We will see in the next theorem that $\sigma^{\D}$ defines an
isomorphism of affine algebraic varieties between the subspaces of
all elements of $k^{\Phi^+}$ satisfying \textup{(R2)} for $\D$
resp. for $\D^{\sigma}$.

By abuse of notation we will denote the images of the reverse root
vectors in the quotient Hopf algebras $u(\D,\mu)$ again by
$x_{ji}$.

\begin{Thm}\label{mainsystem}
Let $\mu \in k^{\Phi^+}$ satisfying \textup{(R2)} for $\D$. Then
\begin{enumerate}
\item $u^{\D^{\sigma}}_{ij}(\sigma^{\D}(\mu)) = (x_{\widetilde{i}\;\widetilde{j}})^N$  for all $1 \leq i < j \leq n+1$.\label{equality}
\item The family $\sigma^{\D}(\mu)$  satisfies \textup{(R2)} for $\D^{\sigma}$ and
$\sigma^{\D^{\sigma}}(\sigma^{\D}(\mu)) =  \mu$.\label{square}
\end{enumerate}
\end{Thm}
\pf (1) Let $1 \leq i<j \leq n+1$. We write
$$\mu'=\sigma^{\D}(\mu).$$
We first compute $u^{\D^{\sigma}}_{ij}(\mu')$. By \eqref{defu} and
\eqref{rule},
\begin{align*}
u^{\D^{\sigma}}_{ij}(\mu')&=\sum_{(i_1,\dots,i_r) \in
I_{ij}}(q-1)^{N(r-2)}\mu'(i_1,\dots,i_r)
(1-(g_{\widetilde{i_r}\;\widetilde{i_{r-1}}})^N).\notag\\
\intertext{By \eqref{defmu'} we can write for all $1 \leq i \leq i_{t-1}<i_{t}\leq j$}
\mu'_{i_{t-1}i_{t}} &= (-1)^{i_{t}-i_{t-1}+1} \tau_{\widetilde{i_t}\;\widetilde{i_{t-1}}}
\sum_{(l^t_1, \dots,l^t_{s_t}) \in I_{i_{t-1},i_t}}(q-1)^{N(s_t-2)}\mu(\widetilde{l^t_{s_t}},\dots,,\widetilde{l^t_1}).\notag\\
\end{align*}
Hence we obtain

\begin{align}
u^{\D^{\sigma}}_{ij}(\mu') &=\sum_{(i_1,\dots,i_r) \in I_{ij}}(q-1)^{N(r-2)} (-1)^{j-i+r-1} \tau(\widetilde{i_r},\dots,\widetilde{i_1})
\notag\\
&\times\left(\sum_{(l^2_1,\dots,l^2_{s_2}) \in I_{i_1i_2}} (q-1)^{N(s_2-2)}\mu(\widetilde{l^2_{s_2}},\dots,\widetilde{l^2_1})\right)\label{left}\\
&\phantom{aa= s_{i_ti_{t+1}}^N}\vdots\notag\\
&\times\left(\sum_{(l^{r}_1,\dots,l^{r}_{s_{r}}) \in I_{i_{r-1}i_r}}
(q-1)^{N(s_{r}-2)}\mu(\widetilde{l^{r}_{s_r}},\dots,\widetilde{l^{r}_1})\right)\notag\\
&\times\left(1-(g_{\widetilde{i_r}\;\widetilde{i_{r-1}}})^N\right).\notag
\end{align}

On the other hand by Theorem \ref{mainreverse}
\begin{align}
(x_{\widetilde{i}\;\widetilde{j}})^N &=  \sum_{(i_1,\dots,i_r) \in I_{ij}} t(\widetilde{i_r}, \dots, \widetilde{i_1}) u_{\widetilde{i_r}\;\widetilde{i_{r-1}}}
\cdots u_{\widetilde{i_2}\;\widetilde{i_1}}.\notag\\
\intertext{By \eqref{defu} and \eqref{widetilde} we have for all $1 \leq i \leq i_{t-1}<i_{t}\leq j$}
u_{\widetilde{i_{t}}\;\widetilde{i_{t-1}}}&=\sum_{(l^t_1, \dots,l^t_{s_t}) \in I_{i_{t-1},i_{t}}}(q-1)^{N(s_t-2)}
\mu(\widetilde{l^t_{s_t}},\dots,\widetilde{l^t_1}) (1-(g_{\widetilde{l^t_2}\;\widetilde{l^t_1}})^N).\notag\\
\intertext{Again we get a large sum of products as before:}
(x_{\widetilde{i}\;\widetilde{j}})^N &=\sum_{(i_1,\dots,i_r) \in I_{ij}} t(\widetilde{i_r}, \dots, \widetilde{i_1})\notag\\
&\times\left(\sum_{(l^2_1,\dots,l^2_{s_2}) \in I_{i_1i_2}} (q-1)^{N(s_2-2)}\mu(\widetilde{l^2_{s_2}},\dots,\widetilde{l^2_1})\right)\label{right}\\
&\phantom{aa= s_{i_ti_{t+1}}^N}\vdots\notag\\
&\times\left(\sum_{(l^{r}_1,\dots,l^{r}_{s_{r}}) \in I_{i_{r-1}i_r}} (q-1)^{N(s_{r}-2)}\mu(\widetilde{l^{r}_{s_r}},\dots,\widetilde{l^{r}_1})\right)\notag\\
&\times \left(1-(g_{\widetilde{l^r_2}\;\widetilde{l^r_1}})^N\right) \cdots \left(1-(g_{\widetilde{l^2_2}\;\widetilde{l^2_1}})^N\right).\notag
\end{align}

The point of the proof is to change the order of the summation
indices. We have to sum over all sequences
$$(l^2_1,\dots,l^2_{s_2},\dots,l^r_1,\dots,l^r_{s_r})\in I_{ij}$$
with $l^2_1=i_1,l^2_{s_2}=i_2=l^3_1, \dots, l^r_{s_r}=i_r$, where each sequence has
$$m=s_2+ \cdots + s_r -r+2$$
elements. Equivalently we can start with an arbitrary sequence
$$(l_1, \dots,l_m) \in I_{ij},$$
then take all subsequences $(l_{k_1}, \dots,l_{k_r})$,
$(k_1,\dots,k_r) \in I_{1,m},$ of\newline $(l_1, \dots,l_m)$ and
define $(i_1,\dots,i_r) \in I_{ij}$ by  $i_p=l_{k_p} $ for all
$p$, \newline $ 1 \leq p \leq r$. Thus the right hand side of
\eqref{left} becomes
\begin{align}
&\sum_{(l_1,\dots,l_m) \in I_{ij}}\sum_{(k_1,\dots,k_r) \in I_{1m}}(-1)^{j-i-r+1}(q-1)^{N(m-2)}\label{R1}\\
&\cdot \mu(\widetilde{l_m},\dots,\widetilde{l_1})
\tau(\widetilde{l_{k_r}},\widetilde{l_{k_{r-1}}},\dots,\widetilde{l_{k_1}})
(1-(g_{\widetilde{l_{k_r}}\;\widetilde{l_{k_{r-1}}}})^N),\notag\\
\intertext{and the right hand side of \eqref{right} becomes}
&\sum_{(l_1,\dots,l_m) \in I_{ij}}\sum_{(k_1,\dots,k_r) \in
I_{1m}}(q-1)^{N(m-r)} \mu(\widetilde{l_m},\dots,\widetilde{l_1})
t(\widetilde{l_{k_r}},\widetilde{l_{k_{r-1}}},\dots,\widetilde{l_{k_1}})
\label{R2}\\
&\cdot
(1-(g_{\widetilde{l_{k_{r-1}+1}}\;\widetilde{l_{k_{r-1}}}})^N)
\cdots
(1-(g_{\widetilde{l_{k_1+1}}\;\widetilde{l_{k_1}}})^N).\notag
\end{align}
Both expressions \eqref{R1} and \eqref{R2} can be simplified.
By Lemma \ref{simplify} we can write in \eqref{R1}
\begin{equation*}
\mu(\widetilde{l_m},\dots,\widetilde{l_1})
\tau(\widetilde{l_{k_r}},\widetilde{l_{k_{r-1}}},\dots,\widetilde{l_{k_1}})=
\mu(\widetilde{l_m},\dots,\widetilde{l_1})\tau_{\widetilde{j}\;\widetilde{i}}.
\end{equation*}
Similarly by Lemma \ref{simplify} and since $N$ is odd we have in \eqref{R2}
\begin{align*}
&\mu(\widetilde{l_m},\dots,\widetilde{l_1})t(\widetilde{l_{k_r}},\widetilde{l_{k_{r-1}}},
\dots,\widetilde{l_{k_1}})\\
&=\mu(\widetilde{l_m},\dots,\widetilde{l_1})(-1)^{j-i-r+1} (q-1)^{N(r-2)} \tau(\widetilde{l_{k_r}},\widetilde{l_{k_{r-1}}},\dots,\widetilde{l_{k_1}})^{-\frac{N-1}{2}} \tau_{\widetilde{j}\;\widetilde{i}}^{\frac{N+1}{2}}\\
&=\mu(\widetilde{l_m},\dots,\widetilde{l_1})(-1)^{j-i-r+1} (q-1)^{N(r-2)} \tau_{\widetilde{j}\;\widetilde{i}}^{-\frac{N-1}{2}} \tau_{\widetilde{j}\;\widetilde{i}}^{\frac{N+1}{2}}\\
&=\mu(\widetilde{l_m},\dots,\widetilde{l_1})(-1)^{j-i-r+1} (q-1)^{N(r-2)} \tau_{\widetilde{j}\;\widetilde{i}}.
\end{align*}

After this simplification we finally obtain
\begin{align}
&u^{\D^{\sigma}}_{ij}(\mu') =
(-1)^{j-i+1} \tau_{\widetilde{j}\;\widetilde{i}}\sum_{(l_1,\dots,l_m) \in I_{ij}}
(q-1)^{N(m-2)}\mu(\widetilde{l_m},\dots,\widetilde{l_1})\label{final1}\\
&\cdot \sum_{(k_1,\dots,k_r) \in I_{1m}}(-1)^r(1-(g_{\widetilde{l_{k_r}}\;\widetilde{l_{k_{r-1}}}})^N),\notag
\intertext{and}
& (x_{\widetilde{i}\;\widetilde{j}})^N =(-1)^{j-i+1} \tau_{\widetilde{j}\;\widetilde{i}}
\sum_{(l_1,\dots,l_m) \in I_{ij}}(q-1)^{N(m-2)}
\mu(\widetilde{l_m},\dots,\widetilde{l_1})\label{final1bis}\\
&\cdot \sum_{(k_1,\dots,k_r) \in I_{1m}}(-1)^r(1-(g_{\widetilde{l_{k_{r-1}+1}}\;\widetilde{l_{k_{r-1}}}})^N) \cdots (1-(g_{\widetilde{l_{k_1+1}}\;\widetilde{l_{k_1}}})^N),\notag
\end{align}
and \eqref{equality} follows from Lemma \ref{group} with
$$h_{st} = (g_{\widetilde{l_t}\;\widetilde{l_s}})^N \text{ for all }1 \leq s<t \leq m$$
and for each sequence $(l_1,\dots,l_m)$.

\medskip

(2) To prove that $\mu'$ satisfies \textup{(R2)} for $\D^{\sigma}$
let $\chi_{ij}' \neq1$. Then
$\chi_{\widetilde{j}\;\widetilde{i}}\neq 1$. Hence for all
$(i_1,\dots,i_r) \in I_{ij}$ we have
$\chi_{\widetilde{i_s}\;\widetilde{i_{s-1}}}\neq 1$ for some $1<s
\leq r$. Thus $\mu(\widetilde{i_r},\dots,\widetilde{i_1}) =0$ ,
and $\mu_{ij}'=0$ by \eqref{defmu'}.

The proof of the equality $\sigma^{\D^{\sigma}}(\sigma^{\D}(\mu))
= \mu$ is similar to the proof of (1). Let $1 \leq i<j \leq n+1$.
By definition
$$\sigma^{\D^{\sigma}}_{ij}(\sigma^{\D}(\mu)) = \tau^{\D^{\sigma}}_{\widetilde{j}\;\widetilde{i}} (-1)^{j-i+1} \sum_{(i_1,\dots,i_r) \in I_{ij}}(q-1)^{N(r-2)} \mu'(\widetilde{i_r},\dots,\widetilde{i_1}).$$
For all $1 \leq i \leq i_{t-1}<i_{t}\leq j$ we have
$$\mu'_{\widetilde{i_t}\;\widetilde{i_{t-1}}} = (-1)^{i_{t}-i_{t-1}+1} \tau_{i_{t-1}i_t}
\sum_{(l^t_1, \dots,l^t_{s_t}) \in I_{i_{t-1},i_t}}(q-1)^{N(s_t-2)}\mu(l^t_1,\dots,l^t_{s_t}).$$
As before for $u^{\D^{\sigma}}_{ij}(\mu')$ we now obtain
\begin{align*}
\sigma^{\D^{\sigma}}_{ij}(\sigma^{\D}(\mu))&=\tau^{\D^{\sigma}}_{\widetilde{j}\;\widetilde{i}}\tau_{ij}
\hspace{-9pt}\sum_{(l_1,\dots,l_m) \in
I_{ij}}(q-1)^{N(m-2)}\mu(l_1,\dots,l_m)
\hspace{-8pt}\sum_{(k_1,\dots,k_r) \in I_{1m}}(-1)^r\notag\\
&=\tau^{\D^{\sigma}}_{\widetilde{j}\;\widetilde{i}}\tau_{ij}\mu_{ij}
\text{\phantom{aaaaa}(by \eqref{ef2})}.\notag
\end{align*}
This proves the claim since
\begin{align*}
\tau^{\D^{\sigma}}_{\widetilde{j}\;\widetilde{i}}&=
\prod_{\widetilde{j} \leq \widetilde{l}<\widetilde{k} < \widetilde{i}}(\chi^{\sigma}_{\widetilde{l}}(g^{\sigma}_{\widetilde{k}}))^N
=\prod_{i <k<l \leq j}(\chi_{l-1}(g_{k-1}))^N
=\prod_{i \leq k<l<j} \chi_{l}^N(g_k),
\intertext{hence}
\tau^{\D^{\sigma}}_{\widetilde{j}\;\widetilde{i}} \tau_{ij}
&= \prod_{i \leq k<l<j}\chi_{l}^N(g_k)\chi_{k}^N(g_l)
=\prod_{i \leq k<l<j} \chi_k(g_k)^{N a_{kl}}
=1.
\end{align*}

\epf

\begin{Cor}\label{sigmaiso}
Let $\mu \in k^{\Phi^+}$ satisfying \textup{(R2)}. Then the map
\begin{equation*}
u(\D^{\sigma},\sigma^{\D}(\mu)) \to u(\D,\mu)
\end{equation*}
given by $x_i^{\sigma} \mapsto x_{\sigma(i)},g \mapsto g, 1 \leq i \leq n ,g \in \G,$
is an isomorphism of Hopf algebras.
\end{Cor}
\pf Let $s_i=1$ for all $1 \leq i \leq n$. Then the triple
$(\id_{\G}, \sigma, (s_i))$ is an isomorphism from
$(\D^{\sigma},\sigma^{\D}(\mu))$ to $(\D,\mu)$ by Theorem
\ref{mainsystem} and Lemma \ref{Lemmaiso1}. Hence the claim
follows from Theorem \ref{cite}. \epf

\section{Hopf algebra isomorphisms}\label{sectioniso}

In this section  let
\begin{align*}
\D &= \D(\G,(g_i)_{1 \leq i \leq n}, (\chi_i)_{1 \leq i \leq n}, (a_{ij})_{1 \leq i,j \leq n}),\\
\D' &= \D(\G',(g'_i)_{1 \leq i \leq n}, (\chi'_i)_{1 \leq i \leq n}, (a_{ij})_{1 \leq i,j \leq n})
\end{align*}
be  data with finite abelian groups $\G$ and $\G'$ and the same
Cartan matrix \eqref{standard} of Type $A_n$. As before $\sigma$
denotes the non-trivial diagram automorphism of $(a_{ij})$ given
by \eqref{defsigma}.

For $s=(s_i)_{1 \leq i \leq n} \in k^n$ and $\mu  \in k^{\Phi^+}$ we define
\begin{align}
s^N&=(s_i^N)_{1 \leq i \leq n}\notag\\
s\cdot\mu&=(s_{ij}\mu_{ij})_{1 \leq i<j \leq n+1}.\notag\\
\intertext{Recall that}
s_{ij} &= \prod_{i \leq l <j} s_l \text{ for all }1 \leq i<j \leq n+1.\notag\\
\intertext{Then}
u_{ij}(s\cdot\mu)&=s_{ij} u_{ij}(\mu) \text{ for all }1 \leq i<j \leq n+1\label{umus}
\end{align}
since $s_{ij}= s_{i_1i_2} \cdots s_{i_{r-1}i_r}$ for all $(i_1,\dots,i_r) \in I_{ij}$.

An isomorphism of data of Cartan type from $\D'$ to $\D$ is a group isomorphism $\varphi : \G' \to \G$ satisfying
\begin{equation}\label{II}
\varphi(g'_i)=g_{i}, \chi'_i=\chi_{i}\varphi\text{ for all } 1
\leq i \leq n.
\end{equation}
We write $\varphi : \D' \xrightarrow{\cong} \D$ if $\varphi$ is an isomorphism from $\D'$ to $\D$.

Note that \eqref{II} implies for all $1 \leq i,j \leq n+1$ that
\begin{align}
\chi'_j(g'_i) &= \chi_j(g_i), \text{ for all }1 \leq i,j \leq n+1,\notag\\
\varphi(g'_{ij}) &= g_{ij} \text{ for all }1 \leq i<j \leq n+1.\notag\\
\intertext{Hence for all $\mu' \in k^{\Phi^+}$}
\varphi(u^{\D'}_{ij}(\mu')) &= u^{\D}_{ij}(\mu') \text{ for all }1 \leq i<j \leq n+1.\label{same}
\end{align}
Let $k^{\times} = k \setminus\{0\}$ denote the multiplicative
group of $k$. In part (b) of (II) below, recall the definition of
$\nu^{\D}$ in \eqref{defTheta}.

\begin{Thm}\label{Hopfiso}
Let $\mu$ and $\mu'$ be families of root vector parameters for $\D$ and $\D'$.
Then the Hopf algebra isomorphisms $u(\D',\mu') \to u(\D,\mu)$ are given by
\begin{alignat*}{2}
\textup{(I) }&x_i \mapsto s_ix_i,\, g \mapsto \varphi(g),\, 1 \leq i \leq n,\,  g \in \G',\text{ where}\phantom{aaaaaaaaaaare gi}\\
&\textup{(a) }\varphi : \D' \xrightarrow{\cong} \D,\text{ and }\\
&\textup{(b) }s=(s_i) \in (k^{\times})^n\text{ such that }\mu'=s^N \cdot \mu,\\
\intertext{and}
\textup{(II) }&x_i \mapsto s_ix_{\sigma(i)},\, g \mapsto \varphi(g),\, 1 \leq i \leq n,\,  g \in \G',\text{ where} \\
&\textup{(a) }\varphi : \D' \xrightarrow{\cong} \D^{\sigma},\text{ and }\\
&\textup{(b) }s=(s_i) \in (k^{\times})^n\text{ such that }\mu'=s^N \cdot \nu^{\D^{\sigma}}(\sigma^{\D}(\mu)).\phantom{}
\end{alignat*}
\end{Thm}
\pf By Theorem \ref{cite} the isomorphisms $u(\D',\mu') \to
u(\D,\mu)$ are given by $x_i \mapsto s_i x_{\rho(i)}, g \mapsto
\varphi(g),1 \leq i \leq n, g \in \G'$, where $\varphi : \G' \to
\G$ is an isomorphism of groups, $\rho = \id$ or $\rho =\sigma$,
and $s=(s_i) \in (k^{\times})^n$ such that
\begin{align}
\varphi(g'_i)&=g_{\rho(i)}, \chi'_i=\chi_{\rho(i)}\varphi\text{ for all } 1 \leq i \leq n.\label{II1}\\
\varphi(u^{\D'}_{ij}(\mu')) &= s_{ij}^N \pi(F^{\rho}(x_{ij}^{\rho})^N) \text{ for all }1 \leq i<j \leq n+1.\label{II3}
\end{align}
We have to show that \eqref{II1} and \eqref{II3} are equivalent to
(I)(a) and (b) if $\rho=\id$, and to (II)(a) and (b) if
$\rho=\sigma$.

By definition \eqref{II1} coincides with (I)(a) if $\rho=\id$ and
with (II)(a) if $\rho = \sigma$.

Let $\rho = \id$ and assume (I)(a). Then the left hand side of
\eqref{II3} is $\varphi(u^{\D'}_{ij}(\mu'))= u^{\D}_{ij}(\mu')$ by
\eqref{same}. For the right hand side of \eqref{II3} we obtain
\begin{alignat*}{2}
s_{ij}^N\pi(F^{\rho}(x_{ij}^{\rho})^N)&=s_{ij}^Nu^{\D}_{ij}(\mu)&&\\
&=u^{\D}_{ij}(s^N\cdot\mu)&\quad&\text{(by \eqref{umus})}.\\
\end{alignat*}
Hence (I)(b) and \eqref{II3} are equivalent by the uniqueness in Lemma \ref{normalization}.

Let $\rho=\sigma$ and assume (II)(a). Let $1 \leq i<j \leq n+1$.
Again it follows from \eqref{same} applied to $\D'
\xrightarrow{\cong} \D^{\sigma}$ that $\varphi(u^{\D'}_{ij}(\mu'))
= u^{\D^{\sigma}}_{ij}(\mu')$. We have shown in Lemma
\ref{Lemmaiso1} and Theorem \ref{mainsystem} that
\begin{alignat*}{2}
\pi(F^{\sigma}(x_{ij}^{\sigma})^N)&=u^{\D^{\sigma}}_{ij}(\sigma^{\D}(\mu)).&&\\
\intertext{Hence}
s_{ij}^N \pi(F^{\sigma}(x_{ij}^{\sigma})^N)&=s_{ij}^N u^{\D^{\sigma}}_{ij}(\sigma^{\D}(\mu))
&\quad&\\
&=s_{ij}^Nu^{\D^{\sigma}}_{ij}(\nu^{\D^{\sigma}}(\sigma^{\D}(\mu)))&&\\
&=u^{\D^{\sigma}}_{ij}(s^N \cdot\nu^{\D^{\sigma}}(\sigma^{\D}(\mu)))&&\text{(by \eqref{umus})}.
\end{alignat*}
Again it follows that (II)(b) and \eqref{II3} are equivalent.
\epf

We know from Theorem \ref{mainsystem} that $\nu^{\D}(\mu)$
satisfies \textup{(R2)} for $\D^{\sigma}$. If we assume (II)(a),
then $\nu^{\D}(\mu)$ also satisfies \textup{(R2)} for $\D'$. Hence
the normalization $\nu^{\D^{\sigma}}(\sigma^{\D}(\mu))$ is a
family of root vector parameters for $\D^{\sigma}$ and $\D'$. In
general $\sigma^{\D}(\mu)$ is not a family of root vector
parameters for $\D^{\sigma}$ since \textup{\textup{(R1)}} is not
necessarily satisfied, and we have to pass to the normalization.

For example
$$\sigma^{\D}_{13}(\mu) = - \tau_{n-1,n+1} (\mu_{n-1,n+1} + (q-1) \mu_{n-1,n} \mu_{n,n+1}),$$
and if $({g^{\sigma}}_{13})^N = g_{n-1,n+1}^N=1$, then
$\mu_{n-1,n+1}=0$, but $\mu_{n-1,n} \mu_{n,n+1}$ can be non-zero
if $g_{n-1}^N\neq 1, g_n^N\neq 1$ and $\chi_{n-1}^N=1,\chi_n^N=1$.
As a realization of this situation take $n=2$ and let  $\G \cong
\mathbb{Z}/(N^2) \times \mathbb{Z}/(N)$  with generators $g$ of
order $N^2$ and $h$ of order $N$. Let $\zeta \in k$ be a root of 1
of order $N$ and $q=\zeta^2$. Define $g_1,g_2$ and characters
$\chi_1,\chi_2$ by
$$g_1=gh,g_2=g^{-1}h,\chi_1(g)=\zeta,\chi_1(h)=\zeta,\chi_2(g)=\zeta^{-2},\chi_2(h)=1.$$
Then $\chi_1(g_1) =q=\chi_2(g_2), \chi_1(g_2)\chi_2(g_1) = q^{-1}$, and
$$g_1^N = g^N \neq 1,g_2^N=g^{-N} \neq1, \chi_1^N=\chi_2^N=1, \text{ and }g_{13}^N= g_1^Ng_2^N=1.$$
Note that \textup{\textup{(R1)}} is trivially satisfied if $g_{ij}^N \neq 1$ for all $1 \leq i<j \leq n+1$.

The next corollary follows immediately from Theorem \ref{Hopfiso}.

\begin{Cor}\label{Coriso}
Let $\mu,\mu'$ be families of root vector parameters for $\D$.
Then the following are equivalent:
\begin{enumerate}
\item $u(\D,\mu') \cong u(\D,\mu)$.
\item There is a family $s \in (k^{\times})^n$ such that \\$\mu' =
\begin{cases}
s^N \cdot \mu,& \text{ if } \D \ncong \D^{\sigma},\\
s^N \cdot \mu \text{ or }s^N \cdot\nu^{\D^{\sigma}}(\sigma^{\D}(\mu)),& \text{ if }\D \cong \D^{\sigma}.
\end{cases}$
\end{enumerate}
\epf\end{Cor}

\begin{Cor}\label{infinite}
Suppose there are $1 \leq i<j \leq n+1,j-i\geq 2$, such that
$g_{ij}^N \neq 1$ and $g_l^N \neq 1, \chi_l^N =1$ for all $i \leq
l <j$. Then the number of isomorphism classes of Hopf algebras of
the form $u(\D,\mu)$ is infinite.
\end{Cor}
\pf By our assumption on $\D$ we can consider families of root
vector parameters $\mu,\mu'$ with $\mu_{l,l+1} = \mu'_{l,l+1}=1$
for all $i \leq l<j$, and with arbitrary elements
$\mu_{ij},\mu'_{ij} \in k$. If $u(\D',\mu') \cong u(\D,\mu)$, then
by Corollary \ref{Coriso} for all $i \leq l<j$ we have
\begin{align*}
\mu'_{l,l+1}&= \begin{cases}
s_l^N \mu_{l,l+1}, & \text{ if }\D \ncong \D^{\sigma},\\
s_l^N \mu_{l,l+1} \text{ or } s_l^N \mu_{\sigma(l),\sigma(l)+1},& \text{ if }\D \cong \D^{\sigma},
\end{cases}\\
\intertext{hence $s_l^N=1$. Thus $s_{ij}^N=1$, and again by Corollary \ref{Coriso} it follows that}
\mu'_{ij} &= \begin{cases}
\mu_{ij},& \text{ if } \D \ncong \D^{\sigma},\\
\mu_{ij} \text{ or }\nu_{ij}^{\D^{\sigma}}(\sigma^{\D}(\mu)),& \text{ if }\D \cong \D^{\sigma}.
\end{cases}
\end{align*}
Hence we obtain infinitely many isomorphism classes of Hopf algebras $u(\D,\mu)$.
\epf

Theorem \ref{Hopfiso} gives the following description of the group
of all Hopf algebra automorphisms of $u(\D,\mu)$.

\begin{Cor}
Let $\mu$ be a family of root vector parameters for $\D$. Then
$\textup{Hopfaut}(u(\D,\mu))$ is isomorphic to the subgroup of
$\Aut(\G) \times (k^{\times})^n$ consisting of all pairs
$(\varphi,s), \varphi \in \Aut(\G), s \in (k^{\times})^n,$ where
\begin{align*}
&\varphi : \D \xrightarrow{\cong} \D, \mu = s^N \cdot \mu, \text{ or }\\
&\varphi : \D \xrightarrow{\cong} \D^{\sigma}, \mu = s^N \cdot \nu^{D^{\sigma}}(\sigma^{\D}(\mu)).
\end{align*}
\epf\end{Cor}

We note the following special case.

\begin{Cor}\label{Corauto}
Let $\mu$ be a family of root vector parameters for $\D$. Then the
group of all Hopf algebra automorphisms of $u(\D,\mu)$ is finite
if $\mu_{i,i+1} \neq 0$ for all $1 \leq i \leq n$. \epf\end{Cor}

\end{document}